\documentclass[a4paper,10 pt]{amsart}

\usepackage[latin1]{inputenc} % support for accents and Umlaute in the text.
\usepackage{amsmath, amssymb, amsfonts, amsthm}
\usepackage[all]{xy}

\newcommand{\C}{\ensuremath{\mathbb{C}}}

\newcommand{\Z}{\ensuremath{\mathbb{Z}}}

\newcommand{\F}{\ensuremath{\mathbb{F}}}

\newcommand{\PP}{\ensuremath{\mathbb{P}}}

\newcommand{\sA}{{\mathcal A}}

\newcommand{\sC}{{\mathcal C}}

\newcommand{\sS}{{\mathcal S}}

\newcommand{\sN}{{\mathcal N}}

\newcommand{\Ga}{\Gamma}

\newtheoremstyle{l-it}% name
     {}%      Space above
     {}%      Space below
     {\itshape}%         Body font
     {}%         Indent amount (empty = no indent, \parindent = para indent)
     {\bfseries}% Thm head font
     {$\;$\textmd{---}}%        Punctuation after thm head
     { }%     Space after thm head: " " = normal interword space;
     {}%         Thm head spec (can be left empty, meaning `normal')

\newtheoremstyle{l-up}% name
     {}%      Space above
     {}%      Space below
     {\upshape}%         Body font
     {}%         Indent amount (empty = no indent, \parindent = para indent)
     {\bfseries}% Thm head font
     {$\;$\textmd{---}}%        Punctuation after thm head
     { }%     Space after thm head: " " = normal interword space;
     {}%         Thm head spec (can be left empty, meaning `normal')

     \newtheoremstyle{citing}% name
    {}%      Space above, empty = `usual value'
    {}%      Space below
    {\itshape}% Body font
    {}%         Indent amount (empty = no indent, \parindent = para indent)
    {\bfseries}% Thm head font
    {$\;$\textmd{---}}%        Punctuation after thm head
    {.35em}%     Space after thm head: " " = normal interword space;
    {\thmnote{#3}}% Thm head spec

\theoremstyle{l-it}
\newtheorem{teo}{Theorem}[section]
\newtheorem{df}[teo]{Definition}
\newtheorem{lem}[teo]{Lemma}
\newtheorem{cor}[teo]{Corollary}
\newtheorem{prop}[teo]{Proposition}

\theoremstyle{l-up}
\newtheorem{ex}[teo]{Example}
\newtheorem{quest}[teo]{Question}

\newtheorem{rem}[teo]{Remark}

\theoremstyle{citing}
\newtheorem{custom}{}

\newcommand{\Proof}{{\it Proof. }}
\newcommand{\QED}{\hspace*{\fill}$\square$\vspace{2ex}}

\newcommand{\hin}{"$\Rightarrow$"$\quad$}
\newcommand{\zur}{"$\Leftarrow$"$\quad$}

\renewcommand{\iff}{\Leftrightarrow}

\newcommand{\ra}{\ensuremath{\to}}

\newcommand{\Aut}{\operatorname{Aut}}
\newcommand{\Int}{\operatorname{Int}}
\newcommand{\Hom}{\operatorname{Hom}}
\newcommand{\Sl}{\ensuremath{\mathbf{Sl}}}

\newcommand{\inverse}[1]{{#1}^{-1}}

\newcommand{\restr}[1]{{\raisebox{-0.3\height}{$\mid_{#1}$}}}

\newcommand{\refb}[1]{{\upshape (\ref{#1})}}

\newcommand{\aslope}[1]{\ensuremath{\mathfrak{a}({#1})}}

\title{Double Kodaira fibrations}
\author{Fabrizio Catanese and S\"onke Rollenske}
\begin{document}

\begin{abstract}
The existence of a  Kodaira fibration, i.e., of a fibration of a compact complex
surface
$S$ onto a complex curve $B$ which is a differentiable but not a holomorphic
bundle, forces the geographical slope $ \nu(S) = c_1^2 (S) / c_2 (S)$ to lie in
the interval
$(2,3)$. But up to now all the known examples had slope $ \nu(S) \leq 2 + 1/3$. 
In this paper we consider a special class of surfaces admitting two such Kodaira
fibrations, and we can construct many new examples,  showing in particular that 
there are such fibrations attaining the slope
$ \nu(S) =  2 + 2/3$. We are able to explicitly describe the moduli space of
such class of surfaces, and we show the existence of Kodaira fibrations which
yield rigid surfaces. We observe an interesting connection between the
problem of the slope of Kodaira fibrations and a 'packing' problem for
automorphisms of algebraic curves of genus $\geq 2$.
\end{abstract}
\maketitle
\rightline{\slshape Version \today}
\tableofcontents
\section{Introduction.}

  It is well known that the Euler characteristic $e$ is multiplicative
for fibre bundles and in 1957 Chern, Hirzebruch and Serre (\cite{chs})
  showed that the same holds true for the signature $\sigma$ if
the fundamental group of the base acts trivially on the cohomology of 
the fibre.

In 1967 Kodaira \cite {kod67} constructed examples of fibrations
of a complex  algebraic surface over a  curve where
multiplicativity of the signature does not hold true,
  and in his honour such fibrations  are nowadays called Kodaira Fibrations.

\begin{df} A {\bf Kodaira fibration} is a
fibration $\psi:S\to B$ of a compact complex
   surface over a compact complex curve, which is a
differentiable but not a holomorphic fibre bundle.
We denote by $b$ the genus of the base curve $B$,
and by $g$ the genus of the  fibre $F$.
\end{df}

It is well known (see section 2) that,
if the fibre genus $g$  is  $ \leq 2$, and there are no
singular fibres,
then one has a holomorphic bundle. Likewise the genus $b$ of the base curve 
of a Kodaira fibration has
to be $\geq 2$.

Atiyah and Hirzebruch (\cite{at69}, \cite{hirz69} ) presented variants of
Kodaira' s construction analysing the relation of the monodromy action to
the non multiplicativity of the signature.

Other constructions of Kodaira fibrations have been later given  by
Gonzalez-Diez and Harvey and others  (see  \cite{gh},
\cite{zaal}, \cite{B-D}
    and references therein) in order to obtain fibrations over curves
of small genus with
fixed signature and fixed fibre genus.

     A precise quantitative  measure of the non multiplicativity
   of the signature is given by  the \emph{geographic slope}, i.e., the ratio
$\nu : =  c_1^2 (S )/ c_2(S)  = K^2_S / e(S)$ between the Chern numbers
of the surface: for Kodaira
fibred surfaces it  lies in the
interval $(2,3)$, in view of the well known Arakelov inequality and of the
improvement by Kefeng Liu (\cite{Liu96}) of the Bogomolov
-Miyaoka-Yau inequality $ K^2_S / e(S)  \leq 3$.

The basic problem we approach in this paper is : which are the
slopes of Kodaira fibrations?

   This problem was posed by Claude Le Brun who
raised the question whether the slopes
can be bounded away from 3:  is it true for instance that
    for a Kodaira fibration the slope is at most 2,91?  In fact,
   the examples by Atiyah, Hirzebruch and Kodaira have a slope
at most  $ 2 + 1/3 = 2,33 \dots $ (see\cite{BHPV}, page. 221)
and if one considers Kodaira fibrations obtained
   from a general complete intersection curve
%under the composition of the Torelli map with the Satake embedding
in the moduli space
$\mathfrak M_g$ of curves of
genus $ g \geq 3$, one obtains a smaller slope (around 2,18).

Our main result in this direction is the following

\begin{custom}[Theorem A] There are Kodaira fibrations with slope equal to $
2 + 2/3 = 2,66\dots $.\end{custom}

Our method of construction is a variant of the one used by Kodaira,
and is briefly described as follows: we consider branched coverings
$ S \ra B_1 \times B_2$ branched on a smooth divisor $D \subset  B_1
\times B_2 $
such that the respective projections $ D \ra B_i$ are \'etale (unramified)
for $ i = 1,2$. We denote these by \emph{double \'etale Kodaira Fibrations}.
The advantage of these is that we are able to completely describe their
moduli spaces.

The starting point is  the topological characterization  of
double Kodaira fibred
surfaces (these are the surfaces admitting two different Kodaira fibrations),
derived from \cite{kot}.

\begin{custom}[Proposition \ref{char}] Let $S$ be a complex surface.
The datum of a
    double Kodaira fibration on $S$ is equivalent to the following data:
\begin{enumerate}
\item  Two exact sequences
\[\xymatrix{1\ar[r] &\Pi_{g_i}\ar[r]& \pi_1(S) \ar[r]^{\bar\psi_i}&
\Pi_{b_i}\ar[r]& 1 & i=1,2}\] where $\Pi_g$ denotes the
fundamental group of a compact curve of genus $g$ and where 
$b_i\geq2, g_i \geq 3$,
such that
\item the composition homomorphism
\[\xymatrix{\Pi_{g_1}\ar[r]& \pi_1(S) \ar[r]^{\bar\psi_2}&
\Pi_{b_2}}\] is neither zero nor injective, and
\item  the Euler characteristic of $S$ satisfies
\[e(S)=4(b_1-1)(g_1-1)=4(b_2-1)(g_2-1).\]
\end{enumerate}
\end{custom}

The above characterization plays an important role in
the   explicit description of the
moduli spaces of Kodaira fibrations.

\begin{custom}[Theorem \ref{moduli}]
Double \'etale Kodaira Fibrations form a closed and open subset
in the moduli spaces of surfaces of general type.
\end{custom}

Thus the moduli space of double \'etale Kodaira fibred
surfaces $S$ is   a union of connected components 
of the moduli spaces of surfaces of general type: we
   conjecture these connected components to be irreducible,
and we prove this conjecture in the  special case of 
\emph{standard} double \'etale Kodaira fibred
surfaces, where we have lots of  concrete examples.

Let us explain how double \'etale Kodaira fibrations can be
constructed starting from
curves with automorphisms, and are indeed related to sets of \'etale
morphisms between two fixed curves.

The simple reason for this is that each component of $D$ is an \'etale covering of
each $B_i$, and thus we can take an \'etale cover $\tilde{B}_1 \to B_1$
dominating each of them; then the pullback $S' \to \tilde{B}_1$ of $S \to B_1$
has the property  that
$D'$ is composed of disjoint graphs of \'etale maps 
$\phi_i:\tilde{B}_1\to B_2$.

The philosophy, as the reader may guess, is then: the larger the
cardinality of $\sS=\{\phi_i\}$
compared to the genus of $B_2$, the bigger the slope, and conversely, once
we find such a set  $\sS$ we get (by the so-called tautological
construction, described
in section \ref{tautological}) plenty of corresponding double (\'etale) Kodaira
fibrations.
If by a further pullback we can achieve $B_1=B_2=B$ and $\sS\subset \Aut(B)$
our question concerning the slope of double Kodaira fibrations is
related to the
following question.

\begin{custom}[Question B] Let  $B$ be a compact complex curve of genus $
b \geq 2$, and let $ \sS \subset \Aut (B)$ be a subset such that all the graphs
$ \Ga_s , s \in \sS$ are disjoint in $ B \times B$: which is  the
best upper bound for
$   | \sS| / (b-1)  $ ?
\end{custom}

We achieve $   | \sS| / (b-1) = 3 $, and in this  way we obtain the
slope $ 8/3$.
Conversely, it is interesting to observe that the cited upper bound
for the slope implies
that $   | \sS| / (b-1)  < 8. $

It would be desirable to find  examples with $   | \sS| / (b-1)  > 3
$, for instance examples with $   | \sS| / (b-1)  =4$ would yield
a  slope equal to $ 2,75$. Even
more interesting would be
to find sharper upper bounds for the slope of Kodaira fibrations.

The consideration of double \'etale Kodaira fibrations
related to curves with many automorphisms enables us also to prove the
following interesting

\begin{custom}[Corollary \ref{rigidexample}]
There are double Kodaira fibred surfaces $S$
which are rigid.\end{custom}

The moduli space of some special Kodaira  fibrations were decribed by  Kas
\cite{Ks} and Jost/Yau \cite{j-y83}; here, we prove the following
general

\begin{custom}[Theorem \ref{standardmoduli}]
The subset of  the moduli space corresponding to
\emph{standard} double \'etale Kodaira fibred
surfaces $S$ (those admitting a pullback
branched in a union of graphs of automorphisms),  is a union of connected components which
   are irreducible, and indeed isomorphic to the moduli
space of   pairs $(B, G)$,  where $B$ is a curve of genus $b$ at
least two and
$G$ is a group
   of biholomorphisms of $B$ of a given topological type.
\end{custom}

\section{General set-up.}

\begin{df}\label{defSKF} A \emph{Kodaira fibration} is a smooth
fibration $\psi_1:S\to B_1$ of a surface over a
curve, which is not a holomorphic fibre bundle.

$S$ is called a \emph{double Kodaira fibred surface} if it admits
a \emph{double Kodaira fibration}, i.e.,  a surjective holomorphic map
$ \psi : S \to B_1 \times B_2$ yielding two  Kodaira fibrations.

    \noindent Let $D\subset B_1\times B_2$ be the branch divisor of $\psi$.  If
both  projections $pr_{B_j}\restr{D}:D\to B_j$ are \'etale we call
$ \psi : S \ra B_1 \times B_2$ a
\emph{double \'etale Kodaira Fibration}.
%\emph{double \'etale Kodaira fibred surface}.
\end{df}

\begin{rem}
Observe that if $S$ admits a double \'etale Kodaira fibration, then
$S$ admits two Kodaira fibrations.
Conversely, if $S$ admits two Kodaira fibrations $\psi_i:S\to B_i , i
= 1,2$ then we consider
the product morphism $\psi_1 \times \psi_2 :S \to B_1 \times B_2$,
and its branch locus $D$.
A calculation in local coordinates shows that at a point $P$ of the
ramification divisor $R$,
there are local coordinates $ (x,y)$ such that $\psi_1 \times \psi_2$
is locally given by
$ (x, x + f(x,y)).$ For instance, if $ f(x,y) = y x^2 - 1/4 y^4$ ,
$R$ is singular at $P$ ($y^3= x^2$) and $D$ is singular
at the image point ($x^8 = z^3$ in suitable coordinates).

Note that a surface $S$ could admit  more than two different Kodaira
fibrations in such a way that  some pair of these  (but not all)
yield a double \'etale Kodaira fibration. 
\end{rem}

\begin{rem}
A. Kas remarked in \cite{Ks} that, if $\phi:S\to B$ is a Kodaira
fibration, then the genus of the base
is at least two and the genus $g$ of the fibre is at least three.
\end{rem}
His argument runs as follows: since $\PP^1$ is the only curve of genus zero
we assume that  $ g \geq 1$.
The fibration
induces a period mapping from the universal cover $\tilde B$ of $B$
to the Siegel upper halfspace of genus $g$, a
bounded subset of $\C^n$ ($ n = \frac{1}{2} g (g+1)$).

If $\tilde B$ is $\PP^1$ or $\C$, i.e., $B$ is rational or elliptic,
such a map has to be constant by
compactness, resp.  by Liouville's theorem. By Torelli's theorem  all
the fibres are isomorphic,
which contradicts our assumption. This
settles the question for the base.

    Since the j-invariant must be constant for a smooth elliptic
fibration, the genus of the fibre is at least two.

Consider now a smooth fibration $\phi:S\to B$ with fibres of genus
two. Then every fibre is a hyperelliptic curve and we
get an induced hyperelliptic involution $\tau$ on $S$. The quotient
$S/\tau$ is a $\PP^1$-bundle over $B$ and the
branch divisor $D$ of the map $S\to S/\tau$ intersects every fibre in
exactly six points.  Hence the restriction of
the projection $S/\tau\to B$ to $D$ is \'etale. Let $\mu:\pi_1(B)\to
S_6$ be the corresponding monodromy homomorphism and let
$f:B'\to B$ the \'etale cover associated to the kernel of $\mu$. We
consider now the pullback of $S/\tau$ by $f$. By
construction the monodromy of $f^*D\to B'$ is trivial hence every
component of $f^*D$ comes from a section $B'\to
f^* (S/\tau)$. This implies that $f^*(S/\tau)$ is isomorphic to the
product $B'\times \PP^1$ where $f^*D$ is given by six
constant sections. Since every fibre of $\psi$ is a double cover of
$\PP^1$ branched over the same six points, all the
fibres are isomorphic and $\psi$ is  a holomorphic bundle.

Hence the fibres of a Kodaira fibration have genus $g$ at
least  three.

In the case of a double Kodaira fibration the genus of the fibre is easily seen
to be at least four by  Hurwitz'  formula, since a fibre of  $\psi_i$ is a
branched covering of a curve of genus at least two, namely the base
curve of the
other fibration.

In particular, $S$ cannot contain rational or elliptic curves, since
no such curve is contained in a fibre or admits
a non-constant map to the base curve. Hence $S$ is minimal and we see, using
the superadditivity of Kodaira dimension, that
$S$ is an algebraic surface of general type.

\begin{rem}\label{nonetale} Let $S$ be a surface admitting two different smooth fibrations
$\psi_i:S\to B_i$ where $b_i := $ genus $(B_i ) \geq 2$  and
also  the fibre genus satisfies
$g_i \geq 2$. If e.g. $\psi_1$ is a holomorphic fibre bundle
map, then $S$ has an \'etale covering which is isomorphic to
 a product of curves, $S\to B_1\times B_2$ is \'etale, and
also $\psi_2$ is a holomorphic fibre bundle.
\end{rem}
\Proof  Let $F$ be a fibre of $\psi_1$. Since the genus of $F$ is at
least two, the automorphism group of $F$ is
finite. Hence we can pull back $S$ by an \'etale map $f:\tilde B_1 \to
B_1$ to obtain a trivial bundle resulting in
the diagram
\[\xymatrix{ \tilde B_1\times F \ar[r]^-{\phi} \ar[d] & S\ar[r]^-{\psi}
\ar[d]^{\psi_1} & B_1\times B_2\ar[dl]\\
\tilde B_1 \ar[r]^f & B_1}\] where $\psi$ is induced by $\psi_1$ and
$\psi_2$. By \cite{Cat00}, Rigidity-Lemma 3.8
there exists a map $g:F\to B_2$ such that $\psi\circ\phi=f\times g.$
Take  $x\in F$ and set  $g(x) : =y$,
$S_y : =\psi_2^{-1}(y)$. In the diagram
\[\xymatrix{\tilde B_1\times\{x\} \ar[r]^-{\phi}\ar[dr]_{f} &
S_y\ar[d]^{\psi\restr{S_y}}\\ & B_1\times \{y\}}\]
    $\phi$ and $f$ are \'etale and consequently also $\psi\restr{S_y}$ is
\'etale. Varying $x$ we see that there can be no
ramification points and $\psi$ and $g$ are \'etale. Now any fibre of
$\psi_2$ is an \'etale covering of $B_1$ of fixed
degree, corresponding to a fixed subgroup of $\pi_1(B_1)$. Thus the
fibres are all isomorphic and we have a
holomorphic bundle. 

\QED

We can now give a topological characterization of double Kodaira
fibrations. We denote by $\Pi_g$ the fundamental
group of a complex curve of genus g.

\begin{prop}\label{char}  Let $S$ be a complex surface. The datum of a
    double Kodaira fibration on $S$ is equivalent to the following data:
\begin{enumerate}
\item  Two exact sequences
\[\xymatrix{1\ar[r] &\Pi_{g_i}\ar[r]& \pi_1(S) \ar[r]^{\bar\psi_i}&
\Pi_{b_i}\ar[r]& 1 & i=1,2}\] with $b_i\geq2, g_i \geq 3$,  and such that
\item the composition map
\[\xymatrix{\Pi_{g_1}\ar[r]& \pi_1(S) \ar[r]^{\bar\psi_2}&
\Pi_{b_2}}\] is neither zero nor injective, and
\item  the Euler characteristic of $S$ satisfies
\[e(S)=4(b_1-1)(g_1-1)=4(b_2-1)(g_2-1).\]
\end{enumerate}
\end{prop}
\Proof  Note that a  holomorphic map $f:C'\to C$ between algebraic
curves is \'etale if and only if the induced map
$f_*$ on the fundamental groups is injective. In fact, in this case
there is a covering space $g:D\to C$
corresponding to the subgroup $f_*(\pi_1(C'))$ in $\pi_1(C)$ and by the lifting
theorem we have a map $\tilde
f:C'\to D$ which induces an isomorphism  of the fundamental
groups. Hence $\tilde f$ is of degree one and
$f=g\circ\tilde f$ is also \'etale. We will apply the previous observation to the map in
({\it ii}\/).

Therefore the only if part of our statement follows from  remark \ref{nonetale}. 

Let's
consider the other direction. Using theorem 6.3
of \cite{Cat03} conditions ({\it i}\/) and ({\it iii}\/) guarantee the
existence of two curves $B_i$ of genus $b_i$
and of holomorphic submersions $\psi_i:S\to B_i$ with ${\psi_i}_* =\bar\psi_i$
whose fibres have respective genera $g_1, g_2$.

Condition ({\it ii}\/) implies that the two fibrations are different
and it remains to see that neither of the
$\psi_i$'s can be a holomorphic bundle. But if it were so, by 
remark \ref{nonetale}, then
  $S \to B_1\times B_2$ would be
\'etale and then the map in ({\it ii}\/) would be injective.

\QED

\begin{rem}\label{net}
Double Kodaira fibrations which are not double \'etale
were constructed in \cite{gh} and \cite{zaal}, essentially
with the same method. The map $F:  B \times B \ra Jac (B)$,
$ (x,y) \mapsto x-y$ contracts the diagonal $
\Delta_B
\subset B \times B$ and maps $ B \times B$ to $ Y : = B - B \subset Jac(B)$.
One takes $\Gamma \subset Y$  to be  a general very ample divisor,
and $D \subset \Gamma \times B$ as $ D : = \cup_{x \in \Gamma} F^{-1}(x) $.
The projection of $D$ to $\Gamma$ is \'etale of degree 2, while
the projection of $D$ to $B$ is of degree equal to $ b : = genus (B)$
but is not
\'etale. The pair $D \subset \Gamma \times B$ yields, as we shall explain in a
forthcoming section, a 'logarithmic Kodaira fibration', and from it one can
construct, via the tautological construction, an actual Kodaira fibration.
\end{rem}

We shall be primarily interested in the case of double \'etale
Kodaira fibrations.
Given a holomorphic map
$\phi$ between two curves let us denote by $\Gamma_\phi$  its graph.
\begin{df}
A double \'etale Kodaira Fibration $S\to B_1\times B_2$ is
said to be  \emph{simple} if 
there exist \'etale maps $\phi_1, \dots, \phi_m$ from $B_1$ to $B_2$
   such that
$D=\dot{ \bigcup}_{k=1,\dots ,m} \Gamma_{\phi_k}$; i.e., if each
component of $D$ is the graph
of one of the $\phi_k$'s.

We say that $S$ is \emph{very simple} if $B_1=B_2$ and all the $\phi_k$'s are
automorphisms.
\end{df}

\begin{lem}\label{pullback}
Every double \'etale Kodaira Fibration admits an \'etale pullback
which is simple.
\end{lem}

\Proof Let  $S\to B_1\times B_2$ be a double \'etale Kodaira Fibration.
The branch divisor $D$  is smooth and we can consider the monodromy
map $\mu:\pi_1(B_1,b_1)\to S_{m_1}$ of the \'etale map $p_1:D\to B_1$. Let
$f:B\to B_1$ the (finite) covering associated to the kernel of $\mu_1$. By
construction the monodromy of the pullback $f^*D\to B$ is trivial, hence every
component maps  to $B$ with degree 1 and the corresponding pullback $f^*S$
is a simple Kodaira fibration.

\QED

\begin{rem}\label{error}
Koll\'ar remarked that it is not always possible
to find a very simple \'etale cover
of $S$. But up to now we do not have a concrete example of this
situation.

\end{rem}

This motivates the following

\begin{df}
A double \'etale Kodaira fibration is called standard
if there exist  \'etale Galois covers $B \to B_i, i= 1,2,$ such that the
\'etale pullback 
\[S' : = S \times_{(B_1 \times B_2) }( B \times B),\]
induced by $ B \times B \to B_1 \times B_2$, is very simple.
\end{df}

\section{Invariants of double \'etale Kodaira Fibrations}\label{invariants}
%\section{Invariants of double \'etale Kodaira fibrations}\label{invariants}
In this section we want to calculate some
invariants of a double \'etale Kodaira fibration. First we need to fix some
notation.

Let $S$ be a double \'etale Kodaira fibration as in Definition \ref{defSKF}.
Let $d$ be the degree of $\psi:S\to B_1\times
B_2$, let  $D\subset B_1\times B_2$  be the branch locus of $\psi$ and let
$D_1,\dots, D_m$ be the connected components of $D$. \

By assumption, the composition map
\[D_i\hookrightarrow B_1\times B_2 \overset{pr_j}{\to} B_j\] is \'etale
and we denote by $d_{ij}$ its degree. Then the
degree of $pr_j\restr{D}:D\to B_j$ is $d_j=\sum_{i=1}^m d_{ij}$ and
we get two formulas for the Euler characteristic of $D_i$,
$$e(D_i)=d_{i1}e(B_1)=d_{i2}e(B_2).$$

The canonical divisor of $B_1\times B_2$ is $K_{B_1\times
B_2}=-e(B_1) B_2 - e(B_2)B_1$ and we
calculate
\begin{align*} K_{B_1\times B_2} \cdot D_i&=-e(B_1) B_2 \cdot  D_i -
e(B_2)B_1 \cdot  D_i\\ &=-e(B_1) d_{i1} - e(B_2)d_{i2}=-2e(D_i)
\intertext{so that by adjunction} D_i^2&=deg (K_{D_i})-K_{B_1\times
B_2} \cdot D_i\\ &=-e(D_i)+2e(D_i)=e(D_i).
\end{align*}

    We write
\[\inverse{\psi}(D_i)=\bigcup_{j=1}^{t_i}R_{ij}\] as a union of
disjoint divisors and denote by $n_{ij}$ the degree
of $\psi\restr{R_{ij}}:R_{ij}\to D_i$ and by $r_{ij}$ the
ramification order of $\psi$ along $R_{ij}$. Then
\[ K_S=\psi^*K_{B_1\times B_2}+\sum_{i,j}(r_{ij}-1)R_{ij}\text{ and
}d=\sum_{j=1}^{t_i} n_{ij}r_{ij}. \]

    To summarize the situation we label the arrows in the following
diagram by the degrees of the corresponding maps:
\[\xymatrix{R_{ij}\ar@{^{(}->}[d]\ar[r]^{n_{ij}}
&D_i\ar@{^{(}->}[d]\ar[r]^{d_{i1}} & B_1\\
    S\ar[r]^-{d}&B_1\times B_2\ar[ru]}\]

    We can now calculate some invariants.

    \begin{prop}\label{formulas}
    In the above situation we have the following formulas
    \begin{enumerate}
\item Setting $\beta_i:=\sum_{j=1}^{t_i}n_{ij}(r_{ij}-1)$
    \begin{gather*}
    c_2(S)=d \ c_2(B_1\times B_2 )-\sum_{i=1}^m \beta_i e(D_i)\\
c_1^2(S)=2c_2(S)-\sum_{i=1}^m\sum_{j=1}^{t_i}
\frac{n_{ij}(r_{ij}-1)(r_{ij}+1)}{r_{ij}}e(D_i)
\intertext{thus the signature  is}
\sigma(S)=\frac{1}{3}(c_1^2(S)-2c_2(S))=
   -\frac{1}{3}\sum_{i=1}^m\sum_{j=1}^{t_i}
\frac{n_{ij}(r_{ij}-1)(r_{ij}+1)}{r_{ij}}e(D_i)
    \end{gather*}
    \item
    \begin{enumerate}
    \item If $\psi:S \to B_1\times B_2$ is a Galois covering then
$r_{ij}=r_i$ and
    \[\frac{c_1^2(S)}{c_2(S)}=2+\frac{-\sum_{i=1}^m
\frac{r_i^2-1}{r_i^2}e(D_i)} {e(B_1)e(B_2)-\sum_{i=1}^m
\frac{r_i-1}{r_i}e(D_i)}.\]
    \item If in addition $D$ is composed of graphs of \'etale maps from $B_1$ to
$B_2$, i.e., $S$ is simple, we have
    \[\frac{c_1^2(S)}{c_2(S)}=2+\frac{1-\frac{1}{m}\sum_{i=1}^m \frac{1}{r_i^2}}
{\frac{2g-2}{m}+1-\frac{1}{m}\sum_{i=1}^m \frac{1}{r_i}}\]
where $g$ is the genus of $B_2$.
    \end{enumerate}
     \end{enumerate}
    \end{prop}
    \Proof
    The first formula  can be obtained by calculating the genus of a
fibre $F$ of $S\to B_1$ using the Riemann-Hurwitz
formula and using $c_2(S)=e(S)=e(B_1)e(F)$.

   For  the second  one a
rather tedious calculation of intersection
numbers is needed so that we prefer to cite \cite{Iz03}\footnote{Note
that we have a slightly different  notation.} which gives us
\begin{gather*}
   \begin{split}
     c_1^2(S)&=d \ c_1^2(B_1\times B_2)-\sum_{i=1}^m\left(2 \ b_i \
e(D_i)+\sum_{j=1}^{t_i}\frac{n_{ij}(r_{ij}-1)(r_{ij}+1)}{r_{ij}}D_i^2\right)\\
     &=2d\  e(B_1) e(B_2) -\sum_{i=1}^m 2 b_i \ e(D_i)-\sum_{i=1}^m\sum_{j=1}^{t_i}
\frac{n_{ij}(r_{ij}-1)(r_{ij}+1)}{r_{ij}}e(D_i)\\
&=2c_2(S)-\sum_{i=1}^m\sum_{j=1}^{t_i}\frac{n_{ij}(r_{ij}-1)(r_{ij}+1)}{r_{ij}}e(D_i).
     \end{split}
     \end{gather*}
The formula for the signature is now obvious.

    Let's look at ({\it ii}\/).  If the covering $\psi$ is Galois the
stabilizers of $R_{ik}$ and $R_{il}$ are
conjugate in the covering group and consequently $n_{ij}=n_i$ and
$r_{ij}=r_i$ do not depend on $j$. Hence for
every $i$
    \[d={t_i} n_{i}r_{i} \iff \frac{{t_i} n_{i}}{d}=\frac{1}{r_i}\text{
and }\beta_i={t_i}n_{i}(r_{i}-1)=d-{t_i}n_{i}\]
    Plugging this into the above formulas we get (a):
    \begin{align*}
    \frac{c_1^2(S)}{c_2(S)}-2&= \frac{-\sum_{i=1}^m \frac{t_i
n_i(r_i-1)(r_i+1)}{r_i}e(D_i)} {de(B_1\times
B_2)-\sum_{i=1}^m {\beta_i}e(D_i)}\\
     &= \frac{-\sum_{i=1}^m \frac{r_i t_i
n_i(r_i-1)(r_i+1)}{r_i^2}e(D_i)} {d\left(e(B_1\times
B_2)-\frac{1}{d}\sum_{i=1}^m (d-t_i n_i)e(D_i)\right)}\\
&=\frac{-\sum_{i=1}^m \frac{r_i^2-1}{r_i^2}e(D_i)}
{e(B_1)e(B_2)-\sum_{i=1}^m \frac{r_i-1}{r_i}e(D_i)}
    \end{align*}
     For (b) we further assume the components of $D$ to have all the
same genus as $B_2$, i.e. $e(D_i)=e(B_2)$ for all $i$. Then
     \begin{align*}
     \frac{c_1^2(S)}{c_2(S)}-2&=\frac{-\sum_{i=1}^m
\frac{r_i^2-1}{r_i^2}e(B_2)} {e(B_1)e(B_2)-\sum_{i=1}^m
\frac{r_i-1}{r_i}e(B_2)}\\
     &=\frac{m-\sum_{i=1}^m \frac{1}{r_i^2}}
{-e(B_1)+m-\sum_{i=1}^m \frac{1}{r_i}}.\\
&=\frac{1-\frac{1}{m}\sum_{i=1}^m \frac{1}{r_i^2}}
{\frac{2g-2}{m}+1-\frac{1}{m}\sum_{i=1}^m
\frac{1}{r_i}}
      \end{align*}  \QED

The above formulas will allow us to give some upper bounds  for the
slope.

\section{Tautological construction}\label{tautological}

\begin{df}
A  pair $(S,D)$ consisting of
\begin{enumerate}
\item a smooth fibration $\psi:S \to B$
  with fibre $F$ and
\item  a divisor $D \subset S$ such that
\begin{enumerate}
\item the projection $D \to B $ is \'etale and
\item  the fibration of pointed curves  $(F, F \setminus D )$ is not isotrivial,
\end{enumerate}
\end{enumerate}
is called  a
{\bf log-Kodaira Fibration}.
\end{df}

Our typical example of the above situation will be a divisor
$$  D  \subset B_1 \times B_2 : = S $$ such that the first projection
$D \to B_1$ is \'etale and the second projection
$D \to B_2$ is nowhere constant.

We shall now see that  in order to construct Kodaira Fibrations it 
suffices to construct
log-Kodaira Fibrations.

\begin{prop}
Let $(S ,D )$ be a log-Kodaira Fibration.

Then for all surjections $\rho:\pi_1(F \setminus D )\to G$ there
  exist an \'etale covering $f:B' \to B$ such that $\rho$ is induced by
a surjection $\pi_1(S'\setminus D') = \pi_1(f^*(S\setminus D))\to G$.
In geometric terms:  for every Galois cover of the fibre $\tilde F\to F$,
  ramified exactly over $F \cap D$,
  there is  an \'etale cover $f:B'\to B$ of the base
  and a Galois covering $\tilde {S'} \to S' :
= f^*S $, ramified over $f^*D$, which extends
$\tilde F\to F$, yielding the following diagram.
\[\xymatrix{\tilde F\ar[d]\ar@{^(->}[rr] &&\tilde{S'}\ar[d]\\
F\ar@{^(->}[r] &S &S'\ar[l]}\]
\end{prop}
\proof

Set for convenience $\hat F : =F\setminus D$ and $\hat S : = S\setminus D$
and assume from now on that $D \neq \emptyset$.  The
fundamental group $\F : =\pi_1(\hat F) $ of the punctured
fibre is then a free group on $2g(F)+D  \cdot F-1$ generators. 
Setting $\Gamma : =\pi_1 (\hat S)$ and $\Pi : =\pi_1(B)$, we have
the exact sequence
\[1\to\F \to \Gamma\to \Pi\to 1\]
    associated to the fibre bundle $\hat S\to B$.

    An \'etale base change $f:B'\to B$ corresponds to a finite index
subgroup $\Pi'\hookrightarrow  \Pi$ and yields a sequence
    $1\to \F \to \Gamma' \to \Pi'\to 1$, associated to
$f^*\hat S\to B'$.

Given instead a finite index
subgroup $\tilde\F\hookrightarrow \F$, corresponding to a branched
cover $\tilde F\to F$ ramified over $F\cap D$, we seek for an exact sequence
\begin{equation} 1\to \tilde\F \to \tilde\Gamma \to \tilde\Pi\to
1,\label{sequence}%\tag{$\ast$}
\end{equation} where $\tilde\Gamma$ and $\tilde\Pi$ are finite index
subgroups in $\Gamma$ and $\Pi$ respectively.
It is necessary and sufficient that $\tilde\Gamma$ be contained in
the normalizer of $\tilde\F$ in $\Gamma$ and
$\F\cap\tilde\Gamma=\tilde\F$.

The proof follows then from the following

\begin{lem} If the covering of the fibre
is Galois, associated to  $1\to \tilde\F\to
\F\overset{\rho}{\to}G{\to} 1$ , then there is always a sequence as in
\refb{sequence}.

\end{lem}

\Proof  Since $\F$ is a normal subgroup, $\gamma \in\Gamma$ operates on $\F$
by
conjugation and hence on $\Hom(\F,G)$ by
$\phi\mapsto \gamma(\phi)=\phi\circ \Int_{\gamma^{-1}}$. Let
$\Gamma_\rho$ be the stabilizer of $\rho$ under this
action.      %%%% $\Gamma_\rho$ of finite index.
    For $\gamma\in\Gamma_\rho$ holds $\rho(\gamma x \gamma^{-1
})=\rho(x)$ and in particular $\gamma$ normalizes
$\tilde\F$, the kernel of $\rho$. Let $\Gamma'$ be the subgroup of
$\Gamma$ generated by $\F$ and $\Gamma_\rho$ . We
have the sequence
\[1\to \F\to \Gamma'\to\Pi'\to 1\] Note that since $\F$ is normal in $\Gamma$
we can write every element $\gamma'\in\Gamma'$ as
a product $\gamma'=fg$ where $f\in\F$ and $g\in\Gamma_\rho$.

Consider a tubular neighborhood $N_0$ around a component $D_0$ of the
pullback of $D$ associated to
$\Pi'\hookrightarrow \Pi$ and let $\gamma_0$ be a small loop around
$D_0$ contained in $N_0\cap F$. We consider
$\gamma_0$ also as an element of $\F$ and regard $N_0$ as a small
neighbourhood of the zero section in the normal
bundle $\sN_{D_0/S}$. The fundamental group $\Pi''$ of $D_0$ is a
group with generators $\alpha_i$, $\beta_i$ and
the single relation $\prod_i [\alpha_i,\beta_i]=1$ and we can write
$\Gamma''=\pi_1(N_0\setminus D_0)$ as a central
extension
\[1\to\Z\to \Gamma'' \to \Pi '' \to 1\] where $\prod_i
[\alpha_i,\beta_i]=\gamma_0^k$ and $k={-D_0^2}$.  This
is proven in \cite{Cat06} with the following argument:
pick a point $P\in D_0$ and write
$D_0= (D_0\setminus P ) \cup \Delta_P$ where
$\Delta_P$ is a small disk around $P$. The $S^1$ bundle
(homotopically equivalent to) $N_0\setminus
D_0$ restricted to these two open subsets is
trivial and the $\sC^\infty$ -Cocycle of $N_0$ with regard to this
trivialisation can be given as $z^k$, where $z$
is a local coordinate in $P$ and $k=-c_1(\sN_{D_0/S})=-D_0^2$.
The  fundamental group of $N_0$ is then calculated
using the Seifert-van Kampen theorem.

By a further base change we may assume that $\gamma_0^k$ is in
$\tilde \F$, e.g. by taking a pullback which makes
$k$ divisible by the order  of $G$ (the exponent of $G$ indeed suffices). The
resulting diagram is
\[\xymatrix{1\ar[r] &\Z\ar[r]\ar[d]& \Gamma'' \ar[r]\ar[d] &
\Pi ''\ar[r]\ar[d]& 1\\
    1\ar[r] &\F\ar[r]& \Gamma' \ar[r]& \Pi'\ar[r]& 1 }\]

Defining $\Gamma'''$ as the inverse image of
$  \Pi ''$ in $\Gamma$ we have
    \[\xymatrix{1\ar[r] &\Z\ar[r]\ar[d]& \Gamma'' \ar[r]\ar[d] &
\Pi ''\ar[r]\ar@{=}[d]& 1\\
    1\ar[r] &\F\ar[r]& \Gamma '''\ar[r]& \Pi '' \ar[r]& 1 }\] and may
finally extend $\rho$ to $ \Gamma'''$ as
follows: After choosing arbitrary images
$\rho''(\alpha_i)=\rho''(\beta_i)$ in $G$ and setting
$\rho''(\gamma_0)=\rho(\gamma_0)$, $\rho$ extends to
$\rho'':\Gamma''\to G$ since every such assignment is
compatible with the relations in the group and the action of $\rho$
on $\gamma_0$.

We are now in the situation
\[\xymatrix{ & \Gamma''\ar[dr]\ar[drr]^{\rho''}\\
\Z\ar[ur]\ar[dr] && \Gamma'''\ar@{-->}[r]& G\\ &
\F\ar[ru]\ar[rru]_\rho}\] and writing each $\gamma\in \Gamma'''$ as
a product $\gamma=fg$ with $f\in \F$ and $g\in\Gamma''$ we set
\[\rho(\gamma)=\rho(f)\cdot\rho''(g).\] To see that this is well
defined let $f_1g_1=f_2g_2$. Then
$\inverse{f_2}f_1=g_2\inverse{g_1}$ is an element of $\F$ and of
$\Gamma''$, i.e. a multiple of $\gamma_0$. Hence
applying $\rho$ and $\rho''$ respectively, we get
$\inverse{\rho(f_2)}\rho(f_1)=\rho''(g_2)\inverse{\rho''(g_1)}$
since the two homomorphisms act in the same way on $\gamma_0$. This
implies $\rho(f_1g_1)=\rho(f_2g_2)$.

It remains to check that this defines a homomorphism. We consider two
elements $f_1g_1$, $f_2g_2$ as above. Since
$\Gamma''$ is contained in $\Gamma'$ we can actually assume that the
$g_i$'s can be written as a combination of the
$\alpha_i$'s, $\beta_i$'s  and are contained in $\Gamma_\rho$, hence
they stabilize $\rho$. Now
\begin{align*}
\rho(f_1g_1f_2g_2 )&=\rho(f_1g_1f_2(\inverse{g_1}g_1)g_2 )&\\ &
=\rho(f_1g_1f_2\inverse{g_1})\rho''(g_1g_2)&\\ &=
\rho(f_1)\rho(g_1f_2\inverse{g_1})\rho''(g_1)\rho''(g_2)&\\ &=
\rho(f_1)\rho(f_2)\rho''(g_1)\rho''(g_2)&\\ &=
\rho(f_1)\rho''(g_1)\rho(f_2)\rho''(g_2)&=\rho(f_1g_1)\rho(f_2g_2)
\end{align*} provided we have chosen the images of the $\alpha_i$'s,
$\beta_i$'s in the centralizer of $G$, which we can
do. The desired group $\tilde \Gamma$ is  the kernel of $\rho:\Gamma'''\to G$.

\QED

If $S$ is a double \'etale Kodaira fibration or a product of
curves, one can easily see that also $S'$ is a
double \'etale Kodaira fibration,  provided that the restriction of the second
projection to $D$ is \'etale. Moreover the following holds:

\begin{lem}\label{vs-stand} Assume that we have a curve $B$ of genus
at least two and a subset
 $\sS=\{\phi_1, \dots, \phi_m\}\subset\Aut B$ such that the graphs of these  automorphisms 
are disjoint subsets of $B\times B$. If we construct  a Kodaira fibration applying the
tautological construction to this log-Kodaira fibration, then the resulting surface 
is in fact a standard Kodaira fibration.
\end{lem}
\Proof
Without loss of generality we may assume that $\phi_1=id_B$, i.e., we identify 
the vertical and the horizontal part of the product $B\times B$ via the 
automorphism $\phi_1$. We fix a base point $x_0$ in $B$. It suffices to prove
 the following: for any \'etale Galois covering $B'\to B$ there exists another 
\'etale covering map $f:B''\to B'\to B$ such that the pullback of 
$D : =\Gamma_{\phi_1}\cup\dots\cup \Gamma_{\phi_m}$ under 
the map $f\times f:B''\times B''\to B\times B$ is composed of graphs of automorphisms of $B''$.

The fundamental group $\pi_1(B,x_0)$ can be considered as a subgroup of a Fuchsian group 
wich acts on the upper half plane. Let $\Gamma$ be the maximal Fuchsian group which
 contains $\pi_1(B,x_0)$ as a normal subgroup. Then we have a sequence
\[ 1\to \pi_1(B,x_0)\to \Gamma\to \Aut(B)\to 1.\]
The Galois covering $B'\to B$ corresponds to an inclusion $\pi_1(B', y_0)\subset\pi_1(B, x_0)$
where $y_0$ maps to $x_0$. Consider the Galois covering $B''\to B'\to B$ associated to the
subgroup $\pi(B'',z_0):= \bigcap_{\gamma\in \Gamma}\gamma \pi_1(B', y_0)\inverse\gamma$ 
which is the largest normal subgroup of $\Gamma$ contained in $\pi_1(B', y_0)$. 
It is in fact a finite index subgroup of $\pi_1(B', y_0)$ since $\pi_1(B', y_0)$ 
is of finite index in $\Gamma$. We have  exact sequences
\[\xymatrix{ &&&1&\\
 1 \ar[r] & \pi_1(B,x_0) \ar[r] &\Gamma \ar[r] & \Aut(B) \ar[r]\ar[u] & 1\\
 1 \ar[r] & \pi_(B'',z_0)\ar@{^(->}[u] \ar[r] &\Gamma \ar@{=}[u] \ar[r] &
 G \ar[r]\ar[u]^\alpha & 1\\
&&& Gal(B''\to B)\ar[u] &\\
&&& 1\ar[u] &}\]
where $G$ is a group of automorphisms of $B''$.

Let $d$ be the degree of the covering $f:B''\to B$. Then the degree of the map
 $f\times f: B''\times B''\to B\times B$ is $d^2$ and it suffices to exhibit for any given automorphism $\phi$ of $B$ a set of $d$ automorphisms of $B''$ such that their graphs map $d$ to $1$ to $\Gamma_\phi$ under the map $f\times f$. In order to do so pick $\psi\in G$ such that $\alpha(\psi)=\phi$ which means $f\circ\psi=\phi\circ f$. Then for any $\sigma\in Gal (B''\to B)$ we have
\begin{align*}
(f\times f)(\Gamma_{\sigma\circ\psi})&=(f\times f)(\{(x,y)\in B''\times B''\mid y=
\sigma\circ\psi(x)\})\\
&= \{(f(x), f(\sigma\circ\psi(x)))\mid x\in B''\}\\
&= \{(f(x), f\circ\psi(x))\mid x\in B''\}\\
&= \{(f(x), \phi(f(x)))\mid x\in B''\} = \Gamma_\phi
\end{align*}
and this map has in fact the same degree as $f$.  \QED

The reason why the monodromy problems mentioned in \ref{error} do not occur 
in this case is that the horizontal and the vertical curve in the product 
are in fact identified via $\phi_1$ and therefore, once we fix a basepoint 
on the curve during the tautological contruction, there is no ambiguity in
 the choice of the basepoint on the other curve.

% 
% \begin{rem} If $S$ is a double \'etale Kodaira fibration or a product of
% curves, one can easily see, that also $S'$ is a
% double \'etale Kodaira fibration,  provided the restriction of the second
% projection to $D$ is \'etale.
% \end{rem}

\section{Slope of double \'etale Kodaira Fibrations} Kefeng Liu proved
in \cite{Liu96} that  the slope $\nu$  of a
Kodaira fibration $S$ satisfies
\[ \nu : = \frac{c_1^2(S)}{c_2(S)}<3\] and LeBrun asked whether the better
bound ${c_1^2(S)}<2.91{c_2(S)}$ would hold.

We will now address the question
about what can be said for double \'etale Kodaira
fibrations. Our purpose here is twofold: to find
effective estimates from below for the maximal slope via the
construction of explicit examples and then to see whether one
can prove also an  upper bound for the slope of double Kodaira
fibrations, using their explicit description.

 To separate the numerical considerations
from the geometrical problems we pose the following
    \begin{df}
    Let $B_1, B_2$ be curves of genus at least two. An admissible
configuration for $B_1\times B_2$ is a tuple
    $\sA=(D,d,\{(t_i,\{r_{ij},n_{ij}\}\})$
     consisting of
    \begin{itemize}
    \item a smooth curve $D=D_1\cup\dots\cup D_m\subset B_1\times B_2$
such that each component $D_i$ maps \'etale to
each of the factors,
    \item a positive integer $d$, and positive integers $t_i$, for all $i=1,\dots, m,$
    \item  for all $i=1,\dots, m,$  a $t_i$-tuple 
$\{(r_{ij}, n_{ij})\}_{j=1,\dots,  t_i}$ of pairs of 
 positive integers with $r_{ij}\geq 2$, and such that
    \[d=\sum_{j=1}^{t_i} n_{ij}r_{ij}.\]
    \end{itemize}
    We call the configuration Galois if $r_{ij}$ does not depend upon $j$,
and we then write
$\sA=(D,d,\{(t_i,r_{i},n_{i})\})$.
    If moreover $D$ is made of graphs of \'etale maps
$\phi_k:B_1\to B_2$ (automorphisms if $B_1
\cong B_2$) we call $\sA$ simple (resp.: very
simple).
    Setting $\beta_i:=\sum_{j=1}^{t_i}n_{ij}(r_{ij}-1)$ we define the
abstract slope of $\sA$ by
\[\aslope{\sA}=2+\frac{-\sum_{i=1}^m\sum_{j=1}^{t_i}\frac{n_{ij}(r_{ij}-1)
(r_{ij}+1)}{r_{ij}}e(D_i)}{d \ e(B_1\times
B_2)-\sum_{i=1}^m {\beta_i} \ e(D_i)}.\]
    \end{df}
    We have seen in section \ref{invariants} that a double \'etale Kodaira
fibration  $S$ gives rise to an admissible
configuration $\sA(S)$. If $\sA$ is any admissible configuration and
$S$ a double \'etale Kodaira fibration with
$\sA(S)=\sA$ we say that $S$ realizes $\sA$. In this case the
abstract slope $\aslope{\sA}$ coincides with the
slope of $S$ by Proposition \ref{formulas}. Note that we also
calculated formulas for the abstract slope of (very) simple Galois
configurations.

    To attain a bound from above for the slope we can now study
independently what is the maximal possible abstract
slope for an admissible configuration and how to realize a given
configuration. We already addressed the second
problem in section \ref{tautological} and proceed by analysing the case of very
simple configurations.

\subsection{Packings of graphs of automorphisms}\label{autos}
In this section we let $B$ be a curve of
genus $g$ and $G=\Aut(B)$ its automorphism group.We want to study
sub\emph{sets} of $G$ such that the corresponding graphs do not intersect. We
can  translate this into a group-theoretical condition:

\begin{lem}\label{groupconditions} Let $P_1,\dots,P_n$ be the points
in $B$ which have a non trivial stabilizer
$\Sigma_{P_i}\subset G$. Let $\pi_i:G \to G/{\Sigma_{P_i}}$ be the
map that sends $\phi\in G$ to the left coset
$\phi\Sigma_{P_i}$.
\begin{enumerate}
\item Two automorphisms $\phi\neq \phi'\in G$ have intersecting
graphs if and only if $\pi_i(\phi)=\pi_i(\phi')$ for
some $i \in \{1,\dots,n\}$.
\item A subset $\sS\subset G$ of cardinality $m$ has nonintersecting
graphs if and only if for each $i \in \{1,\dots,n\}$the image of $\sS$ under
the map
\[\pi_i : G\to G/{\Sigma_{P_i}}\qquad
g\mapsto g\Sigma_{P_i}\] has
cardinality $m$. In particular:
\[m\leq Minimum\{|G/ \Sigma_{P_i}|\}_{i=1,\dots,n}\]
\end{enumerate}
\end{lem} Note that, if $Q_1,\dots,Q_r$ are the branch points of the
quotient map $B\to B/G$ and $P_i'$ ($\forall i=1,\dots r$) is an arbitrary
point in the inverse image of $Q_i$ , then the non trivial
stabilizers of points are exactly all the
subgroups conjugated to the stabilizers $\Sigma_{P_i'}$ $(i=1,\dots,r)$.

\Proof Let $\phi,\phi'$ be two automorphisms of $B$. Their graphs
intersect in some point $(P,Q)\in B\times B$ iff
$\phi(P)=\phi'(P)=Q$. But this means $\phi^{-1}\circ \phi'(P)=P$, i.e.,
$\phi^{-1}\circ \phi'\in \Sigma_P$ or
equivalently $\phi\Sigma_{P}=\phi'\Sigma_{P}$.

\QED

It is now a natural question to ask for the maximal possible $m$ that one can 
realize,  given a curve $B$, or given a fixed genus $b$ (of $B$). 

For the formulation of a partial
result we introduce the following notation: we say that $B$ is of
type $(\nu_1,\dots,\nu_k)$ if $B/G$ has genus zero
and the map  $B\to B/G$ is a ramified covering, branched over $k$
points with respective multiplicities $\nu_i$. 

We always order the branch points so that $\nu_1\leq\dots\leq\nu_k$.
\begin{prop}\label{table}
\begin{enumerate}
\item If the genus $g$ of $B$ is at least two, the maximal cardinality
  $m$ of a subset with nonintersecting graphs is smaller or equal to
$3(g-1)$ unless the type of $B$ occurs in the following table:\\
\begin{minipage}{\textwidth}
\begin{center}
\begin{tabular}{l|c|c} type & upper bound for $m$& $|G|$\\
\hline (2,2,2,3)& $4(g-1)$ & $12(g-1)$\\ (2,3,7) & $12(g-1)$ &
$84(g-1)$\\ (2,3,8) & $6(g-1)$ & $48(g-1)$\\ (2,3,9)
& $4(g-1)$ & $36(g-1)$\\ (2,4,5) & $8(g-1)$ & $40(g-1)$\\ (2,4,6) &
$4(g-1)$ & $24(g-1)$\\ (2,5,5) & $4(g-1)$ &
$20(g-1)$\\ (3,3,4) & $6(g-1)$ & $24(g-1)$\\
\end{tabular}
\end{center}
\end{minipage}
\item If the genus of $B$ is small we get the following list:

\begin{minipage}{\textwidth}
\begin{center}
\begin{tabular}{l|c|c} type & upper bound for $m$& up to genus\\
    \hline (2,2,2,3)& $2(g-1)$ & $30$\\
     (2,3,7) & $3(g-1)$ & $23$\\
(2,3,8) & $3(g-1)$ & $23$\\
   (2,3,9) & $2(g-1)$ & $23$\\
    (2,4,5) & $2(g-1)$ & $23$\\
     (2,4,6) & $2(g-1)$ & $50$\\
      (2,5,5)&$4/3(g-1)$ & $50$\\
      (3,3,4) & $3(g-1)$ & $50$\\
\end{tabular}
\end{center}
\end{minipage}
\end{enumerate}
\end{prop} If the genus of the curve is one, we can clearly produce
such an arbitrarily large  subset by choosing
appropriate translations.

\Proof
    Part (i) is a case by case analysis using the previous Lemma. Let $B$
be a curve of genus $g\geq 2$, let $G$ be its
automorphism group and  let $h$ be the genus of $B/G$. Let $P_1,\dots,
P_k\in B/G$ be the branch points and
$\nu_1\leq\dots\leq	\nu_k$ be the corresponding ramification
indices. Then we have the Hurwitz formula
\[2g-2=|G|\left(2h-2+\sum_{i=1}^{k}(1-\frac{1}{\nu_i})\right)\] and by
the lemma a maximal subset as above has at most
cardinality
\[\mu:=\frac{|G|}{\nu_k}=\frac{2g-2}{\nu_k\left(2h-2+\sum_{i=1}^{k}(1-\frac{1}{\nu
_i})\right)},\]
where we set
$\nu_1=1$ if there is no ramification. Note that the denominator can
never be zero since this would imply $g=1$. We
distiguish the following cases:
\begin{itemize}
\item[$h\geq 2$:] Clearly
\[\mu\leq
\frac{2g-2}{\nu_k\left(2+\sum_{i=1}^{k}(1-\frac{1}{\nu_i})\right)}\leq
g-1.\]
\item[$h=1$:] We have \[\mu\leq
\frac{2g-2}{\nu_k\sum_{i=1}^{k}(1-\frac{1}{\nu_i})} \leq 2(g-1)\]
\item[$h=0$:] Also in this case we necessarily have ramification and
\[\mu\leq \frac{2g-2}{\nu_k\left(-2+
\sum_{i=1}^{k}(1-\frac{1}{\nu_i})\right)},\] hence we have to check in
which
cases holds
\[0<\lambda:=\nu_k\left(-2+
\sum_{i=1}^{k}(1-\frac{1}{\nu_i})\right)<\frac{2}{3}\] Since $k\geq 5$
implies
$\lambda\geq 1$ we have $k$ at most 4, and 
$\nu_k > 2$. If $k=4$ then
$\lambda\geq (1/2) \nu_k - 1$ thus $\nu_k = 3$ and one sees immediately
that $(2,2,2,3)$ is the
only possibility. If $k=3$ one can check that $1-
\sum_{i=1}^{3}{1}/{\nu_i}\geq 1-{1}/{2}-{1}/{3}-{1}/{7}={1}/{42}$,
(which corresponds to $|G|=84(g-1)$), hence there are only finitely
many cases for $\nu_k$ which are easy to consider and which
 yield exactly the remaining cases in
the above table.
\end{itemize}

For part ({\it ii}\/) note that a finite group $G$ can occur as an
automorphism group of a curve of type
$(\nu_1,\dots,\nu_k)$ iff there are distinct elements $g_1, \dots,
g_k$  in $G$ such that $g_1, \dots,
g_{k-1}$ generate $G$, $\prod_{i=1}^{k}g_i=1$ and the order of $g_j$
is $\nu_j$. (cf. section
\ref{bound} for a construction.)  For all possible combinations of
groups and generators up to the given genus,
maximal subsets satisfying the conditions of the above Lemma were
calculated using the program GAP (cf. \cite{gap}).

\QED

\begin{rem}\label{[24,3]} The bounds in the second table are sharp,
that is, there exist examples that realize the
given upper bound. The smallest group realizing $3(g-1)$ is
$\Sl(2,\Z/3\Z)$ acting on a curve of genus 2 of type $(3,3,4)$.

We will see in Remark \ref{betterbound} that the slope inequality
obtained by Liu implies in fact the better bound $m<8(g-1)$.
\end{rem}
\begin{quest}\label{packingquestion} It is clear that we can realize the bound
$m=3(g-1)$
for arbitrary large genera $g$ by taking Galois \'etale coverings of the
 examples we have obtained.

 Can one prove that $3(g-1)$ is  an upper bound
for all curves? 

%The maximal possible $m$
%seem only to depend on the abstract automorphism group, not on the  particular
%generators chosen. Could it be helpful to phrase
%the problem in terms of Fuchsian groups?
\end{quest}

\subsection{Bounds for the slope}\label{bound}
Since the slope of a Kodaira fibration does not change under \'etale
pullback, by lemma \ref{pullback} it suffices to treat  the slope for a simple
configuration.  We do this here for the Galois case.

\begin{prop}\label{slopebound} Let $\sA=(D_1\cup\dots\cup D_m, d,
\{t_i, r_i, n_i\})$ be a simple, Galois configuration  and let $g$ be the genus
of the target curve $B_2$. If $m\leq 3(g-1)$ then
$\aslope{\sA}\leq 2 + 2/3$ with equality if and
only if $m=3(g-1)$ and all the ramification indices $r_i$ are equal to three.
\end{prop}

Note that  we do not know any example of a possible (very) simple
configuration with $m>3(g-1)$.  

We  believe that the same result as above should hold 
also in the non Galois case.

\Proof First of all let's assume that $m=3(g-1)$ and let us calculate
$\aslope{\sA}-8/3$ in this case.
\begin{align*}
\aslope{\sA}-8/3&=\frac{1-\frac{1}{m}\sum_{i=1}^m \frac{1}{r_i^2}}
{\frac{2}{3}+1-\frac{1}{m}\sum_{i=1}^m
\frac{1}{r_i}}-\frac{2}{3}\\ &=\frac{3-\frac{3}{m}\sum_{i=1}^m
\frac{1}{r_i^2}-\frac{10}{3}+\frac{2}{m}\sum_{i=1}^m
\frac{1}{r_i}} {5-\frac{3}{m}\sum_{i=1}^m \frac{1}{r_i}}
\intertext {and if we denote by $m_k$ the number of components $D_i$
of $D$ which have ramification index $r_i=k$}
&= \frac{-\frac{1}{3}+\frac{1}{m}\sum_{k}\left(
\frac{2m_k}{k}-\frac{3m_k}{k^2}\right)} {5-\frac{3}{m}\sum_{k}
\frac{m_k}{k}}=\frac{-\frac{1}{3}+\frac{1}{m}\sum_{k}m_k\frac{2k-3}{k^2}}
{5-\frac{3}{m}\sum_{k} \frac{m_k}{k}}
\end{align*} The expression $\frac{2k-3}{k^2}$ has a global maximum
in $k=3$ and hence for the numerator
\[-\frac{1}{3}+\frac{1}{m}\sum_{k}m_k\frac{2k-3}{k^2}\leq-\frac{1}{3}+\frac{1}{m}
\sum_{k}m_k\frac{3}{9}=0\]

with
equality if and only if $m_3=m$ and all other  $m_k$ 's are zero. Consequently
$\aslope{\sA}-8/3\leq 0$ with equality if and only if all
ramification indices are three. It remains to show that
the abstract slope can only decrease if $m < 3(g-1)$ which follows by
induction from the next lemma.

\QED

\begin{lem} Let $\sA=(D_1\cup\dots\cup D_{m+1}, d, \{t_i, r_i,
n_i\})$ be a simple, Galois configuration with
$m\leq 4(g(B_2)-1)$ and let $\sA'=(D_1\cup\dots\cup D_m, d, \{t_i, r_i,
n_i\})$ be the configuration obtained by omitting the
last component. Then $\aslope{\sA'}<\aslope{\sA}$.
\end{lem}

\proof Using again the formulas from proposition \ref{formulas} we calculate
\[\aslope{\sA'}-2<\aslope{\sA}-2\]
\begin{multline*}
\iff \left(m-\sum_{i=1}^{m}\frac{1}{r_i^2}\right)\left(2g-2+m+1
-\sum_{i=1}^{m+1}\frac{1}{r_i}\right)<\\
\left(m+1-\sum_{i=1}^{m+1}\frac{1}{r_i^2}\right)\left(2g-2+m
-\sum_{i=1}^{m}\frac{1}{r_i}\right)
\end{multline*}
\begin{gather*}
\iff \left(m-\sum_{i=1}^{m}\frac{1}{r_i^2}\right)\left(1-\frac{1}{
r_{m+1}}\right)<
\left(1-\frac{1}{r_{m+1}^2}\right)\left(2g-2+m
-\sum_{i=1}^{m}\frac{1}{r_i}\right)\\
\iff \aslope{\sA'}-2=\frac{1-\frac{1}{m}\sum_{i=1}^m \frac{1}{r_i^2}}
{\frac{2g-2}{m}+1-\frac{1}{m}\sum_{i=1}^m
\frac{1}{r_i}}<\frac{1-\frac{1}{r^2_{m+1}}}{1-
\frac{1}{r_{m+1}}}=1+\frac{1}{r_{m+1}}
\end{gather*} The denominator on the left  is bigger or equal to one
since $\frac{2g-2}{m}\geq \frac{1}{2}$ and
$r_i\geq 2$. Hence the left hand side is smaller than one which is
strictly smaller than the right hand side and we
are done.\QED

\begin{ex}\label{example}
We want now to construct an example of a double Kodaira fibration
which actually realizes the slope $8/3$  thereby proving Theorem A.
First of all we construct the curve mentioned in  Remark \ref{[24,3]}.

Let $P_1$, $P_2$, $P_3$ be distinct points in $\PP^1$ and let $\gamma_1$,
$\gamma_2$, $\gamma_3$ be simple geometrical loops around
these points. The fundamental group  $\pi_1(\PP^1\setminus\{P_1, P_2,
P_3\})$ is generated by the $\gamma_i$'s  with
the relation $\gamma_1\gamma_2\gamma_3=1$. Consider in
$\Sl(2,\Z/3\Z)$ the elements
\[g_1=\begin{pmatrix}0 & 2\\1 &2\end{pmatrix},\quad
g_2=\begin{pmatrix}0 &1\\2 &2\end{pmatrix},
\quad g_3=\begin{pmatrix}2 & 2\\2 &1\end{pmatrix}\]
    and define  $\rho:\pi_1(\PP^1\setminus\{P_1, P_2, P_3\})\to
\Sl(2,\Z/3\Z)$ by $\gamma_i\mapsto g_i$.
     This map is well defined and surjective, because $g_1$ and $g_2$
generate $\Sl(2,\Z/3\Z)$ and $g_1g_2g_3=1$. We
define $B$ to be the ramified Galois cover of $\PP^1$ associated to
the kernel of $\rho$. By construction
$\Sl(2,\Z/3\Z)$ acts on  $B$ as the Galois group of the covering and
by the Riemann-Hurwitz formula
\[g(B)=\frac{|\Sl(2,\Z/3\Z)|}{2}\left(\sum_{i=1}^3\left(1-\frac{1}
{ord(g_i)}\right)-2\right)+1=\frac{24}{2}\left(1-\frac{1}{3}-\frac{1}
{3}-\frac{1}{4}\right)+1=2\]
The subset
\[\sS=\left\{ \phi_1=\begin{pmatrix}1 &0\\0 &1\end{pmatrix},
\phi_2=\begin{pmatrix}2 & 0\\1 &2\end{pmatrix},
\phi_3=\begin{pmatrix}0& 1\\2 &1\end{pmatrix} \right\} \subset
\Sl(2,\Z/3\Z)\] satisfies the conditions of  Lemma
\ref{groupconditions} since $\phi_2$,  $\phi_3$ and $\phi_3 \circ
\inverse{\phi}_2$ have no fixed points being of order six and hence
gives us $3=3(g(B)-1)$ graphs of  automorphisms in $B \times B$ which do not
intersect. We denote the corresponding divisor by $D$.

In order to use the tautological construction we  have to construct a
ramified covering of a curve of
genus two minus three points (which  we denote for the sake of simplicity by $(B\setminus D)$)
and Proposition \ref{slopebound} tells
us that the ramification indices should all be equal to 
three.

Let $\alpha_1,\beta_1, \alpha_2, \beta_2$ be generators for
$\pi_1(B)$ and let $\gamma_1$, $\gamma_2$, $\gamma_3$ simple geometrical
loops around the three points such that
\[\pi_1(B\setminus D)=<\alpha_1,\beta_1, \alpha_2, \beta_2,
\gamma_1,\gamma_2,\gamma_3>/({\Pi[\alpha_i,\beta_i]=\gamma_1\gamma_2\gamma_3})\]

is a free group and we can define a
map
\begin{align*}
\rho:\pi_1(B\setminus D)&\to \Z/3\Z\\
\gamma_i&\mapsto 1\\
\alpha_i, \beta_i&\mapsto 0
\end{align*} which induces the desired ramified covering $F\to B$.

At this point we  use the tautological construction,
but we observe that in this case only the first
\'etale covering $B' \to B$ is needed.

Indeed, the divisor $D = D_1 + D_2 + D_3$ has degree 3 on each fibre
of the first projection $p : B \times B \to B$ and the homomorphism
$\rho$ determines a simple cyclic covering of the fixed fibre
$B^0 : = \{x_0\} \times B$, ramified on the divisor $D \cap B^0$.

Therefore there is a divisor $M$ on $B \cong B^0$ such that 
the simple cyclic covering is obtained by taking the cubic root of $D$
in the line bundle corresponding to $M$, and in particular the
following linear equivalence holds:
$$ 3 M \equiv D | _{B^0}.$$
This linear equivalence determines $M$ up to 3-torsion,
and the monodromy of $M$ is the same as the monodromy of $\rho$.

Therefore, if we take as before the \'etale covering $ B' \to B$
associated to the stabilizer of $\rho$, and denote by $D'$ the pull
back of $D$ on $ B' \times B$, then on $ B'
\times B$ the divisor
$D' - 3 p_2^* (M)$ is trivial on each fibre of the first projection
$p_1$,whence there is a divisor $L'$ on $B'$ such that
$D - 3 p_2^* (M) = p_1^* (L')$.

By intersecting with the fibres of the second projection we find that
$ deg (L') = 0$, hence there is a divisor $M'$ on $B'$ such that $ L'
\equiv 3 M'$, and we conclude that on  $ B' \times B$ we have the linear
equivalence
$$D' \equiv  3  ( p_2^* (M) +  p_1^* (L') )$$
and we can take the corresponding simple cyclic covering branched on
$D'$ 
 inside the line bundle
corresponding to the divisor $  p_2^* (M) +  p_1^* (L')$.

We obtain in this way a  double \'etale Kodaira fibration which is in fact 
a Standard Kodaira Fibration by Lemma \ref{vs-stand}. In particular we have 
a Kodaira fibration with base curve $B'$  and with fibre of genus $g=7$
(since $2 g = 2 =  3 \cdot 2 + 3 \cdot 2  $).

Since the associated ramified covering is branched
exactly over $D'$ with ramification index three at each component,
the formula for the slope of a simple configuration calculated in
Proposition \ref{formulas} yields

\[\frac{c_1^2(S)}{c_2(S)}=2+\frac{1-\frac{1}{3}\sum_{i=1}^3
\frac{1}{3^2}} {-\frac{e(B)}{3}+1-\frac{1}{3}\sum_{i=1}^3
\frac{1}{3}}=\frac{8}{3}.\]
\end{ex}

\begin{rem}\label{betterbound}
We can also use this construction to give a partial answer to the question
raised in \ref{packingquestion}. Knowing that the slope of a Kodaira fibration
is strictly smaller than 3 it follows that $m<8(g-1)$. In fact, via a
suitable base change we obtain a divisor 
$D' \subset B' \times B$ such that
\begin{enumerate}
\item  if $m$ is odd, then there is a component $D_1$ mapping to $B'$ with
degree one,
\item setting $D'' := D'$ if $m$  is even, and $D'' := D ' - D_1$ if $m$  is
odd, then
\item we can take a double cover branched over $D''$.
\end{enumerate}
The Kodaira fibration constructed in this way turns out, under the 
assumption
$m\geq 8(g-1)$, and in view of the above  formulas, to have a slope 
$\geq 3$: this is a contradiction.

It follows in particular as a consequence: if $B$ is a curve of genus $2$ and
we have 8 \'etale maps from a fixed curve $C$ of arbitrary genus to $B$,
then two of
them have a coincidence point.
\end{rem}

\section{The moduli space}

This section is devoted to the description of the moduli space of
double \'etale
Kodaira Fibrations. We start with some lemmas.

\begin{lem}\label{lem1}
Let $B_1$, $B_2$ be curves of genus $b_i\geq 2$ resp. and let $C\subset B_1\times
B_2$ be an irreducible curve. Then  
\begin{itemize}
\item
$C$ is smooth and the restricted projections
$p_i: C\to B_i$ are \'etale if and only if 
\item
the negative of
the selfintersection of $C$ attains its maximum possible value, i.e., iff
\[-C^2=2m_i(b_i-1)\qquad  (i=1,2)\]
where $m_1=C \cdot \{*\}\times B_2$ and $m_2=C \cdot B_1\times\{*\}$.
\end{itemize}
\end{lem}
\Proof
\hin We calculated this at the begining of section \ref{formulas}.\\
\zur Let $p=p(C)$ be the arithmetic genus of $C$. Then
\begin{align*}
2p-2=& K_{B_1\times B_2} \cdot C+C^2
=2(b_1-1)m_1+2(b_2-1)m_2 - 2(b_j-1)m_j\\=&2m_i(b_i-1)\qquad (i\neq j)
\end{align*}
Let $\tilde C  \to C$ be the normalisation and let $g=g(\tilde C)$ be the geometric
genus of $C$. We have $2p-2\geq 2g-2$ by the normalization sequence and on
the other hand $2g-2\geq 2m_i(b_i-1)=2p-2$ by the Hurwitz formula. Hence $g=p$,
C is smooth and equality holds in the last inequality, i.e., there is no
ramification and the maps $p_i$ are \'etale.

\QED
\begin{rem}
In general we see that  $K_{B_1\times B_2} \cdot C+C^2=2m_i(b_i-1) + 
2\delta+\rho_i$
where $\delta$ is the 'number of double points' and $\rho_i$ is the total
ramification index of $ C \to B_i$. So
\[-C^2= 2m_j(b_j-1)-2\delta-\rho_i\qquad (i\neq j)\]
\end{rem}

\begin{lem}
Assume that we have a family of effective divisors $(D_t)_{t\in T}$,
$ D_t \subset (B_{1,t} \times B_{2,t})$, such that the
special fibre $D: =D_0=nC$ with $C$ as in Lemma \ref{lem1}. If $D' $ is
another fibre ($D' = D_t $ for some $t$), then
$D'$  is of the same type $D'=nC'$ (the integer $n$ being the same as before).
\end{lem}
\Proof
Write $D'=\sum_j r_j C_j$ as a sum of irreducible components, so that
$C_i \cdot C_j\geq 0$ for $i\neq j$. Write also
$m_1^j=C_j \cdot\{*\}\times B_{2,t}$ and
$m_2^j=C_j \cdot B_{1,t}\times\{*\}$. We calculate
\begin{align*}
-{D'}^2&=\sum_j r_j^2(-C_j^2)-2\sum_{i\neq j}r_i r_j C_i \cdot C_j\\
& \leq \sum_j r_j^2 (-C_j^2)\leq \sum_j r_j^2 2m_i^j(b_i-1)
\intertext{and also}
-{D'}^2 &= n^2(-C^2)=n^2 2m_i(b_i-1).
\end{align*}
Hence the following conditions hold:
\[\sum_j r_j^2 m_i^j\geq n^2 m_i\qquad, nm_i = \sum_j r_j m_i^j.\]
Since $D$ is the special fibre every component $C_i $ tends to a positive
multiple of $C$, we have $m_i^j\geq m_i$ and putting together the two
inequalities yields
\[ \sum_j r_j^2{ m_i^j}^2\geq n^2 m_i^2= \left(\sum_j r_j
m_i^j\right)^2\geq\sum_j r_j^2{ m_i^j}^2\]
therefore in fact equality holds, there is only one summand and $D'=r_0C_0$. To
conclude the proof we look again at the conditions
\[n^2m_i \leq r_0^2m_i^0 , \qquad nm_i=r_0 m_i^0 , \qquad m_i\leq m_i^0.\]
Combining the two inequalities with the equality in the middle  we
get $n\leq r_0\leq n$ and we are done by observing that also $C'=C_0$ fullfills the
conditions of Lemma \ref{lem1}.

\QED

\begin{teo}\label{moduli}
Being a double \'etale Kodaira Fibration is a closed and open
condition in the moduli
space.
\end{teo}

\Proof
Due to the previous Lemma it remains to show the closedness. Assume then that we
have a 1-parameter family of surfaces with general fiber $S_t$ a double \'etale
Kodaira fibration. By the topological characterization (Proposition \ref{char})
also the special fibre $S_0$ is a double Kodaira fibration. Moreover, 
 by Lemma \ref{pullback}, we  may assume that $S_t$ is a branched covering of
$B_{1,t}\times B_{2,t}$ branched over $D_t=\sum_i k_i D_{i,t}$, where the
$D_{i,t}$'s are disjoint graphs of \'etale maps $\phi_i:B_{1,t}\to B_{2,t}$.

Now, $S_0\to B_{1,0}\times B_{2,0}$ is branched over $D_0:=\sum_i k_i
\nu_i D_{i,0}$ where $D_{i,t}$ tends to $\nu_i D_{i,0}$. Since however
$D_{i,t}.(B_{1,t}\times\{*\})=1$ we have $\nu_i D_{i,0}.(B_{1,0}\times\{*\})=1$
which implies $\nu_i=1$.  Hence $D_{i,0}$ is the graph of a map
$\phi'_i:B_{1,0}\to B_{2,0}$ and another application of Lemma \ref{lem1} shows
that also $\phi_i'$ is \'etale and $S_0$ is a double \'etale Kodaira
Fibration.

\QED

We can now describe the moduli space of Standard Kodaira Fibrations in detail.
Let $S$ be a Standard Kodaira Fibration: then there exists  a minimal common
 Galois cover  $B'$ of $B_1, B_2$ yielding   an \'etale pullback $S'$
which is very simple.  We call $B'$ the {\em simplifying covering curve}. We have diagrams
\[\xymatrix{S'\ar[rr]^{\psi_2'}\ar[dr]^\pi \ar[dd]^{\psi_1'} &&
B'\ar[d]^{f_2}\\
%f_1^*S\ar[r]\ar[d]
& S\ar[r]^{\psi_2}\ar[d]^{\psi_1} &B_2 \\
B'\ar[r]^{f_1} & B_1}
\qquad \xymatrix{B'\times B'\ar[r]& B_1\times B_2\\
D'=\bigcup_{\phi\in \sS} \Gamma_\phi \ar@{^{(}->}[u]\ar[r] &D\ar@{^{(}->}[u]}
\]
where $D$ is  the ramification divisor of $\psi=\psi_1\times \psi_2$
and $D'=\pi^*D$ is made of the graphs of a set of automorphisms
$\sS\subset \Aut(B')$. If we denote the Galois group of $f_i$ by
$G_i$  $(i=1,2)$  the following holds:

\begin{teo}\label{standardmoduli} Let $S$ be a standard Kodaira fibred
surface and let $\mathfrak N$ be the
irreducible (and connected) component of the moduli space containing $[S]$.
$\mathfrak N$ is then isomorphic to the moduli space of the pair
$(B', G)$, where $B'$ is the simplifying covering curve  defined above and $G$ is the
subgroup of $\Aut(B')$ generated by $G_1$, $G_2$ and $\sS$.
\end{teo}
\proof
Let us first consider the case where $S=S'$, i.e., where $S$ itself is
very simple. By proposition \ref{char} every  deformation in the 
large of $S'$ is a
branched cover af a product surface $ B_1 \times B_2$. Moreover, clearly
$ B_1 = B_2$ if (*) there is a component of the branch locus mapping to both
curves $B_1, B_2$ with degree $1$. So let
$S_t, t \in T$, be a family with connected parameter space $T$, having
$S'$ as a fibre. It is clear that the set of points of $T$ where (*) 
holds is open.
It is also closed because in the proof
of theorem
\ref{moduli} we have seen that the type  of the branch divisor remains the same
under specialization.

  We have seen that $\mathfrak N$ parametrizes surfaces
which are  very simple and indeed a branched covering
of a product $ B \times B$ branched on the union of graphs of automorphisms.

   The automorphisms defining the components of the branch divisors in
different fibers are clearly pairwise isotopic to each other and
therefore we obtain a family of curves with automorphisms.

For each curve let $G$ be the finite group generated by these
automorphisms.  This group has a faithful representation on the
fundamental group of the curve, and
therefore the group $G$ remains actually constant.

  $G$ is a finite
group and we have a faithful  action on Teichm\"uller space 
$\mathfrak T_b$. We use
now Lemma 4.12 of
\cite{Cat00} (page 29) to the effect that the fixed locus of this action is
a connected submanifold (diffeomorphic to an Euclidean space), hence
the moduli space of such pairs (B,G) is irreducible.

Viceversa any element in this moduli space gives rise to a complex
structure on the differentiable manifold underlying $S'$.

Consider now the general case. It is clear that any
deformation of $S$ induces a deformation of $B' \to B_i$, hence any deformation
of $S$ yields a deformation of the pair $  (B', G)$.

Conversely, any deformation of the pair $ (B', G)$
yields a deformation of the pair $ D' \subset B' \times B' $
such that the group $ G_1 \times G_2$ leaves $D'$  and the monodromy 
of the unramified covering of $ ( B' \times B') - D' $ invariant.

 \QED

\begin{cor}\label{rigidexample}
There exist  double \'etale Kodaira fibred surfaces which  are rigid.

\end{cor}
\Proof
Take the fibration constructed in Example
\ref{example}: the automorphisms corresponding to the ramification 
divisor generate
the whole triangle group of type (3,3,4) and it is well known that
pairs $(B,G)$ yielding a triangle curve are rigid.
Similarly for the other examples in
proposition \ref{table} which yield $ m = 3 (g-1)$.

\QED

\section*{Acknowledgements}

The present research took place in the realm of the D.F.G. Schwerpunkt
" Globalse Methoden in der komplexen Geometrie".

The first author would like to thank Joe Harris for an interesting conversation,
Janos Koll\'ar for remark
\ref{error}, and Gianpietro Pirola for remark \ref{net}.

The second author wishes to thank all the members of the Lehrstuhl Mathematik
II at the University of Bayreuth, especially Ralf Gugisch,  for useful
discussions concerning the computational problems.

\vfill

\noindent
{\bf Author's address:}

\bigskip

\noindent
Dipl. Math.  S\"onke Rollenske \\
Lehrstuhl Mathematik VIII\\
    Mathematisches Institut\\
Universit\"at Bayreuth, NWII\\
     D-95440 Bayreuth, Germany

e-mail: Soenke.Rollenske@uni-bayreuth.de

\bigskip
\noindent
Prof. Dr. Fabrizio Catanese\\
Lehrstuhl Mathematik VIII\\
   Mathematisches Institut\\
Universit\"at Bayreuth, NWII\\
     D-95440 Bayreuth, Germany

e-mail: Fabrizio.Catanese@uni-bayreuth.de

\end{document}